December 15, 2009

# SERIES TRANSFORMATION FORMULAS OF EULER TYPE, HADAMARD PRODUCT OF FUNCTIONS, AND HARMONIC NUMBER IDENTITIES

**Khristo N. Boyadzhiev**
Ohio Northern University
Department of Mathematics
Ada, Ohio 45810, USA
k-boyadzhiev@onu.edu

**Abstract.** The integral representation of the Hadamard product of two functions is used to prove several Euler-type series transformation formulas. As applications we obtain three binomial identities involving harmonic numbers and an identity for the Laguerre polynomials. We also evaluate in a closed form certain power series with harmonic numbers.

1. **Introduction**

Given a function

$$f(t) = \sum_{k=0}^{\infty} a_k t^k \tag{1.1}$$

holomorphic in a neighborhood of zero, Euler's series transformation formula says that

$$\frac{1}{1-t} f\left(\frac{t}{1-t}\right) = \sum_{n=0}^{\infty} t^n \left\{ \sum_{k=0}^{n} \binom{n}{k} a_k \right\} \tag{1.2}$$

(see [1] or [5]). This formula can be used, among other things, to evaluate the binomial expression

$$c_n = \sum_{k=0}^{n} \binom{n}{k} a_k \quad (n = 0,1,2...), \tag{1.3}$$

which is called the *binomial transform* of the sequence $\{a_k\}$. Some examples can be found in [1]. With the substitution

$$\frac{t}{1-t} = z, \tag{1.4}$$

formula (1.2) takes the form

$$f(z) = \frac{1}{z+1} \sum_{n=0}^{\infty} \left(\frac{z}{z+1}\right)^n \left\{\sum_{k=0}^{n} \binom{n}{k} a_k\right\}. \tag{1.5}$$

We want to point out here a useful connection between this transformation formula and the integral representation of the Hadamard product for series (see [11, section 4.6])

$$\sum_{n=0}^{\infty} a_n b_n z^n = \frac{1}{2\pi i} \oint_L g(\frac{z}{\lambda}) f(\lambda) \frac{d\lambda}{\lambda}, \tag{1.6}$$

where $L$ is an appropriate closed curve around the origin and

$$g(t) = \sum_{k=0}^{\infty} b_k t^k \tag{1.7}$$

is a second analytic function defined on a neighborhood of zero. By choosing the function $g(t)$ appropriately, one can generate transformation formulas like (1.5). In this paper we prove several such formulas and present a number of applications. Our main results are propositions 1, 3, 4, 5 and corollaries 1 and 2, all in section 2. Corollary 2, for instance, gives the evaluation in a closed form of the harmonic number sums

$$\sum_{n=0}^{\infty} H_n n^m z^n \tag{1.8}$$

and more generally, of the sums,

$$\sum_{n=0}^{\infty} (H_{p+n} - H_p) \binom{p+n}{n} n^m z^n, \tag{1.9}$$

for any integer $m \geq 0$.

In the other three sections we present several more applications – identities for the harmonic numbers and for the Laguerre polynomials.

Let $L_n(x)$ be the Laguerre polynomials. In section 3 we prove the identity

$$\int_0^x \frac{L_n(t)-1}{t} dt = \sum_{k=1}^{n} \frac{L_k(x)-1}{k} \quad . \tag{1.10}$$

for $n = 1, 2, ...,$ .

We shall use the extended harmonic numbers

$$H_\alpha = \psi(\alpha+1) + \gamma, \ \text{Re}\, \alpha > -1, \tag{1.11}$$

where $\psi(z)$ is the digamma (or psi) function and $\gamma = -\psi(1)$ is Euler's constant [10]. When $\alpha = k \geq 0$ is an integer, then

$$H_k = 1 + \frac{1}{2} + ... + \frac{1}{k}, \ H_0 = 0 \tag{1.12}$$

are the usual harmonic numbers. In section 4 we prove a symmetric identity involving these numbers. Namely, for any $n \geq 0$ and $\text{Re}\, \alpha > -1$, we have

$$\sum_{k=0}^{n} \binom{n}{k}\binom{\alpha}{k} H_k = \binom{\alpha+n}{n}(H_\alpha + H_n - H_{\alpha+n}) \ . \tag{1.13}$$

In section 5 two more identities are proved:

$$\sum_{k=0}^{n} \binom{\alpha+1}{n-k}\binom{\alpha+k}{k}(-1)^{n-k} H_k = \frac{1}{n} + \binom{\alpha}{n}\frac{(-1)^{n-1}}{n}, \tag{1.14}$$

$$\sum_{k=0}^{n} \binom{\alpha}{n-k}\binom{\alpha+k}{k}(-1)^{n-k} H_k = H_n + \sum_{k=1}^{n} \binom{\alpha}{k}\frac{(-1)^{k-1}}{k} \ . \tag{1.15}$$

## 2. Euler-type series transformations

Throughout the paper, the two functions $f(t),\ g(t)$ are defined by the series (1.1) and (1.7) above. We shall need an important lemma.

**Lemma 1.** The following representation holds

$$\sum_{n=0}^{\infty} b_n\, h(z)^n \left\{ \sum_{k=0}^{n} \binom{n}{k} a_k \right\} = \frac{1}{2\pi i} \oint_L g(h(z)(1+\frac{1}{\lambda})) f(\lambda) \frac{d\lambda}{\lambda} \ , \tag{2.1}$$

where $h(z)$ is an appropriate function for which the above expression is defined and the integral is a Cauchy type integral on a closed curve around the origin as in (1.6).

*Proof.* By Cauchy's integral formula, for $k = 0, 1, ...,$

$$a_k = \frac{1}{2\pi i} \oint_L \frac{f(\lambda)d\lambda}{\lambda^{k+1}} . \tag{2.2}$$

Then obviously,

$$\sum_{k=0}^{n} \binom{n}{k} a_k = \frac{1}{2\pi i} \oint_L \left(1+\frac{1}{\lambda}\right)^n f(\lambda) \frac{d\lambda}{\lambda}. \tag{2.3}$$

Multiplying this identity by $b_n h(z)^n$ and summing for $n$ we obtain (2.1).

We shall compare now (2.1) to (1.6). Our purpose is, by choosing $g(t)$ and $h(z)$ appropriately, to make the integrand in (2.1) look like the integrand in (1.6). This means the function $g(h(z)(1+\frac{1}{\lambda}))$ should be written in terms of function(s) in $z$ and function(s) in $z/\lambda$ only. This technique is demonstrated in Proposition 1.

**Proposition 1.** Let $\alpha$ be any complex number. Then the following representation extending (1.5) is true

$$\sum_{n=0}^{\infty} \binom{\alpha}{n}(-1)^n a_n z^n = (z+1)^\alpha \sum_{n=0}^{\infty} \left(\frac{z}{z+1}\right)^n \binom{\alpha}{n}(-1)^n \left\{\sum_{k=0}^{n} \binom{n}{k} a_k\right\} . \tag{2.4}$$

*Proof.* Choose

$$h(z) = \frac{z}{z+1}, \; g(t) = (1-t)^\alpha = \sum_{n=0}^{\infty} \binom{\alpha}{n}(-1)^n t^n, \; b_n = \binom{\alpha}{n}(-1)^n . \tag{2.5}$$

A simple computation shows that

$$g\left(h(z)\left(1+\frac{1}{\lambda}\right)\right) = (z+1)^{-\alpha}\left(1-\frac{z}{\lambda}\right)^\alpha , \tag{2.6}$$

and (2.1) takes the form

$$\sum_{n=0}^{\infty} b_n \left(\frac{z}{z+1}\right)^n \left\{\sum_{k=0}^{n} \binom{n}{k} a_k\right\} = \frac{(1+z)^{-\alpha}}{2\pi i} \oint_L \left(1-\frac{z}{\lambda}\right)^\alpha f(\lambda)\frac{d\lambda}{\lambda}. \tag{2.7}$$

This representation in view of (1.6) leads to (2.4).

When $\alpha = -1$ we have

$$\binom{-1}{n} = (-1)^n \tag{2.8}$$

and (2.4) becomes (1.5). The proof is completed.

Formula (2.4) originates from Euler and can be found in [8, p. 294].

**Proposition 2.** The following (exponential) version of Euler's series transformation formula holds

$$\sum_{n=0}^{\infty} \frac{a_n}{n!} z^n = e^{-z} \sum_{n=0}^{\infty} \frac{z^n}{n!} \left\{ \sum_{k=0}^{n} \binom{n}{k} a_k \right\}. \tag{2.9}$$

*Proof.* Take $h(z) = z$, $g(t) = e^t$. Then (2.1) becomes

$$\sum_{n=0}^{\infty} \frac{z^n}{n!} \left\{ \sum_{k=0}^{n} \binom{n}{k} a_k \right\} = \frac{e^z}{2\pi i} \oint_L e^{\frac{z}{\lambda}} f(\lambda) \frac{d\lambda}{\lambda}, \tag{2.10}$$

and (2.9) follows from (1.6).

This exponential transformation formula can be found in [5] with a simple direct proof.

In the next two examples we use the natural logarithmic function. To our knowledge, these results are new. In all expansions we assume that $|z|$ is small enough to secure convergence.

**Proposition 3.** With $f(t)$ as in (1.1) the following representation holds

$$f(0)\log(1+z) + \sum_{n=1}^{\infty} \frac{z^n}{n} a_n = \sum_{n=1}^{\infty} \left( \frac{z}{z+1} \right)^n \frac{1}{n} \left\{ \sum_{k=0}^{n} \binom{n}{k} a_k \right\}. \tag{2.11}$$

*Proof.* In formula (2.1) we put

$$h(z) = \frac{z}{z+1}, \ g(t) = -\log(1-t) = \sum_{n=1}^{\infty} \frac{t^n}{n}, \ b_n = \frac{1}{n}. \tag{2.12}$$

Then

$$g\left( h(z)\left( 1 + \frac{1}{\lambda} \right) \right) = \log(1+z) - \log\left( 1 - \frac{z}{\lambda} \right), \tag{2.13}$$

and the right hand side in (2.1) becomes

$$\log(1+z)\frac{1}{2\pi i}\oint_L f(\lambda)\frac{d\lambda}{\lambda} - \frac{1}{2\pi i}\oint_L \log(1-\frac{z}{\lambda})f(\lambda)\frac{d\lambda}{\lambda} . \qquad (2.14)$$

Obviously, (2.11) follows from here. The first integral in (2.14) is $a_0 \log(1+z)$. The proof is completed.

We now use the series expansion [7, (7.43), p. 351-352]

$$g(t) = \frac{-\log(1-t)}{(1-t)^{p+1}} = \sum_{n=0}^{\infty}(H_{p+n} - H_p)\binom{p+n}{n}t^n, \qquad (2.15)$$

where $\operatorname{Re} p > -1$ and

$$H_p = \psi(p+1) + \gamma \qquad (2.16)$$

are the extended harmonic numbers - see (1.11).

**Proposition 4**. For any $p$ with $\operatorname{Re} p > -1$ we have

$$\sum_{n=0}^{\infty}(H_{p+n} - H_p)\binom{p+n}{n}a_n z^n + \log(1+z)\sum_{n=0}^{\infty}\binom{p+n}{n}a_n z^n$$

$$= \frac{1}{(1+z)^{p+1}}\sum_{n=0}^{\infty}\left(\frac{z}{z+1}\right)^n (H_{p+n} - H_p)\binom{p+n}{n}\left\{\sum_{k=0}^{n}\binom{n}{k}a_k\right\} . \qquad (2.17)$$

*Proof.* With $g(t)$ as above in (2.15), and

$$h(z) = \frac{z}{z+1} , \qquad (2.18)$$

a simple computation gives

$$g\left(h(z)\left(1+\frac{1}{\lambda}\right)\right) = (1+z)^{p+1}\left[\frac{\log(1+z)}{(1-\frac{z}{\lambda})^{p+1}} - \frac{\log\left(1-\frac{z}{\lambda}\right)}{(1-\frac{z}{\lambda})^{p+1}}\right], \qquad (2.19)$$

and therefore, the right hand side of (2.1) becomes

$$(1+z)^{p+1} \left[ \frac{\log(1+z)}{2\pi i} \oint_L \frac{1}{(1-\frac{z}{\lambda})^{p+1}} \frac{f(\lambda)d\lambda}{\lambda} - \frac{1}{2\pi i} \oint_L \frac{\log(1-\frac{z}{\lambda})^{p+1}}{(1-\frac{z}{\lambda})^{p+1}} \frac{f(\lambda)d\lambda}{\lambda} \right]. \quad (2.20)$$

According to Hadamard's formula (1.6) this equals the left hand side of equation (2.17) multiplied by the factor $(1+z)^{p+1}$. Note that

$$\frac{1}{(1-t)^{p+1}} = \sum_{n=0}^{\infty} \binom{p+n}{n} t^n. \quad (2.21)$$

The entire identity (2.17) now follows from Lemma 1, i.e. from (2.1), as in this case

$$b_n = (H_{p+n} - H_p)\binom{p+n}{n}. \quad (2.22)$$

The proof is completed.

When $p = 0$, then $b_n = H_n$ and we have the corollary:

**Corollary 1.** With $f(t)$ as in (1), the following series transformation formula holds

$$\sum_{n=0}^{\infty} H_n a_n z^n + \log(1+z) f(z) = \frac{1}{1+z} \sum_{n=0}^{\infty} \left(\frac{z}{z+1}\right)^n H_n \left\{\sum_{k=0}^{n} \binom{n}{k} a_k\right\}. \quad (2.23)$$

In our next corollary we present one interesting particular case of (2.17). Let $m \geq 0$ be an integer. Taking $a_k = (-1)^k k^m$ and changing $z$ to $-z$ in (2.17) we obtain the identity

$$\sum_{n=0}^{\infty} (H_{p+n} - H_p)\binom{p+n}{n} n^m z^n + \log(1-z) \sum_{n=0}^{\infty} \binom{p+n}{n} n^m z^n$$

$$= \frac{1}{(1-z)^{p+1}} \sum_{n=0}^{\infty} \left(\frac{z}{1-z}\right)^n (H_{p+n} - H_p)\binom{p+n}{n} (-1)^n \left\{\sum_{k=0}^{n} \binom{n}{k}(-1)^k k^m\right\}. \quad (2.24)$$

We shall use now the representation

$$\sum_{n=0}^{\infty} \binom{p+n}{n} n^m z^n = \frac{1}{(1-z)^{p+1}} \omega_{m,p+1}\left(\frac{z}{1-z}\right) \quad (2.25)$$

from [2]. Here $\omega_{m,p+1}$ are the generalized geometric polynomials

$$\omega_{m,p+1}(x) = \sum_{k=0}^{m} S(m,k)\Gamma(p+k+1)x^k, \tag{2.26}$$

introduced in [2], and $S(m,k)$ are the Stirling numbers of the second kind [3], [7]. Also, in the right hand side of (2.24) we shall use the well known representation of the numbers $S(m,k)$

$$(-1)^n n! S(m,n) = \sum_{k=0}^{n} \binom{n}{k}(-1)^k k^m, \tag{2.27}$$

to obtain the following.

**Corollary 2.** For every $\operatorname{Re} p > -1$ and for every integer $m \geq 0$, with the polynomials $\omega_{m,p+1}$ defined in (2.26), we have the following harmonic number summation

$$\sum_{n=0}^{\infty}(H_{p+n} - H_p)\binom{p+n}{n}n^m z^n \tag{2.28}$$

$$= \frac{1}{(1-z)^{p+1}}\left\{-\log(1-z)\,\omega_m\!\left(\frac{z}{1-z}\right) + \sum_{n=0}^{m}\left(\frac{z}{1-z}\right)^n (H_{p+n} - H_p)\binom{p+n}{n}n! S(m,n)\right\}.$$

Notice that the sum on the right hand side is now finite, as $S(m,n) = 0$ for $n > m$. When $p = 0$, (2.28) turns into the shorter summation formula

$$\sum_{n=0}^{\infty} H_n n^m z^n = \frac{1}{1-z}\left\{-\log(1-z)\,\omega_m\!\left(\frac{z}{1-z}\right) + \sum_{n=0}^{m}\left(\frac{z}{1-z}\right)^n H_n n! S(m,n)\right\}. \tag{2.29}$$

Here $\omega_m = \omega_{m,1}$, that is,

$$\omega_m(x) = \sum_{k=0}^{m} S(m,k) k! x^k. \tag{2.30}$$

The polynomials $\omega_m$ are the geometric polynomials as defined in [2].

The representations (2.28) and (2.29) can also be derived directly from the summation formula (3.23) in [2]. This was done by Dil and Kurt in the recent paper [4].

For completeness, at the end of this section we obtain from (1.6) another interesting series transformation formula, also attributed to Euler.

**Proposition 5.** Given two analytic functions $f(t)$ and $g(t)$ as in (1.1) and (1.7), the following representation is true.

$$\sum_{n=0}^{\infty} a_n b_n t^n = \sum_{n=0}^{\infty} \frac{g^{(n)}(-t)}{n!} t^n \left\{ \sum_{k=0}^{n} \binom{n}{k} a_k \right\}. \tag{2.31}$$

This transformation formula can be found in a modified form in [9, Chapter 6, problem 19, p. 245].

*Proof.* Multiplying both sides in equation (2.3) by

$$\frac{g^{(n)}(-t)}{n!} t^n \tag{2.32}$$

and summing for $n$ we obtain

$$\sum_{n=0}^{\infty} \frac{g^{(n)}(-t)}{n!} t^n \left\{ \sum_{k=0}^{n} \binom{n}{k} a_k \right\} = \frac{1}{2\pi i} \oint_L \sum_{n=0}^{\infty} \frac{g^{(n)}(-t)}{n!} \left( \frac{t}{\lambda} + t \right)^n f(\lambda) \frac{d\lambda}{\lambda} \tag{2.33}$$

$$= \frac{1}{2\pi i} \oint_L g\left(\frac{t}{\lambda}\right) f(\lambda) \frac{d\lambda}{\lambda},$$

by recognizing the Taylor expansion of $g\left(\frac{t}{\lambda}\right)$ centered at $-t$ inside the first integral. Now (2.31) follows from (1.6) and the proof is completed.

Note that (2.31) can also be written in the form

$$\sum_{n=0}^{\infty} a_n b_n t^n = \sum_{n=0}^{\infty} \frac{(-1)^n g^{(n)}(t)}{n!} t^n \left\{ \sum_{k=0}^{n} \binom{n}{k} (-1)^k a_k \right\}. \tag{2.34}$$

3. **An identity for the Laguerre polynomials**

To illustrate how the above formulas work we shall provide some examples. The first application involves the Laguerre polynomials $L_n(x) = \frac{e^x}{n!} \frac{d^n}{dx^n} (e^{-x} x^n)$, see [10].

**Corollary 3.** Let $L_n(x)$ be the Laguerre polynomials. Then for $n = 1, 2, ...,$ we have

$$\int_0^x \frac{L_n(t) - 1}{t} dt = \sum_{k=1}^{n} \frac{L_k(x) - 1}{k}. \tag{3.1}$$

*Proof.* Take an arbitrary $x$ and set in (2.11)

$$f(t) = \sum_{k=0}^{\infty} \frac{(-xt)^k}{k!} = e^{-xt}, \quad a_k = \frac{(-x)^k}{k!}. \tag{3.2}$$

Thus

$$\log(1+z) + \sum_{n=1}^{\infty} \frac{z^n (-x)^n}{n! n} = \sum_{n=1}^{\infty} \left(\frac{z}{z+1}\right)^n \frac{1}{n} L_n(x), \tag{3.4}$$

because of the well-known identity [10, p. 153],

$$\sum_{k=0}^{n} \binom{n}{k} \frac{(-x)^k}{k!} = L_n(x). \tag{3.5}$$

With the substitution (1.4) equation (3.4) can be written in the form

$$-\log(1-t) + \sum_{n=1}^{\infty} \left(\frac{t}{1-t}\right)^n \frac{(-x)^n}{n! n} = \sum_{n=1}^{\infty} \frac{t^n}{n} L_n(x), \tag{3.6}$$

where we divide both sides by $1-t$ to obtain

$$\frac{-\log(1-t)}{1-t} + \frac{1}{1-t} \sum_{n=1}^{\infty} \left(\frac{t}{1-t}\right)^n \frac{(-x)^n}{n! n} = \sum_{n=1}^{\infty} t^n \sum_{k=1}^{n} \frac{L_k(x)}{k}, \tag{3.7}$$

(for the last equality we use property (5.5) below). From (3.5), dividing by $x$ and integrating we find

$$\sum_{k=1}^{n} \binom{n}{k} \frac{(-x)^k}{k! k} = \int_0^x \frac{L_n(t) - 1}{t} dt, \tag{3.8}$$

and therefore, from (1.2)

$$\frac{1}{1-t} \sum_{n=1}^{\infty} \left(\frac{t}{1-t}\right)^n \frac{(-x)^n}{n! n} = \sum_{n=1}^{\infty} t^n \int_0^x \frac{L_n(t) - 1}{t} dt. \tag{3.9}$$

At the same time

$$\frac{-\log(1-t)}{1-t} = \sum_{n=1}^{\infty} H_n t^n. \tag{3.10}$$

Substituting (3.9) and (3.10) in (3.7) and comparing coefficients we arrive at (3.1).

## 4. A symmetric identity for harmonic numbers

First we obtain an equivalent version of the representation (2.15).

**Corollary 4.** For any complex $\alpha$ and $|t| < 1$,

$$\frac{-\log(1-t)}{(1-t)^{\alpha+1}} = \sum_{n=0}^{\infty} t^n \left\{ \binom{\alpha+n}{n} H_n - \sum_{k=1}^{n} \binom{n}{k}\binom{\alpha}{k} H_k \right\}. \tag{4.1}$$

*Proof.* Set in (2.23)

$$f(z) = (1+z)^\alpha = \sum_{n=0}^{\infty} \binom{\alpha}{n} z^n, \quad a_k = \binom{\alpha}{k} \tag{4.2}$$

to get

$$\sum_{n=0}^{\infty} \binom{\alpha}{n} H_n z^n + (1+z)^\alpha \log(1+z) = \frac{1}{1+z} \sum_{n=0}^{\infty} \left(\frac{z}{z+1}\right)^n H_n \binom{\alpha+n}{n}, \tag{4.3}$$

since, by the well-known Vandermonde identity [6],

$$\sum_{k=0}^{n} \binom{n}{k}\binom{\alpha}{k} = \binom{\alpha+n}{n}. \tag{4.4}$$

With the substitution (1.4) we can write (4.3) in the form

$$\frac{1}{1-t} \sum_{n=0}^{\infty} \left(\frac{t}{1-t}\right)^n \binom{\alpha}{n} H_n - \frac{\log(1-t)}{(1-t)^{\alpha+1}} = \sum_{n=0}^{\infty} H_n \binom{\alpha+n}{n} t^n. \tag{4.5}$$

Next, for the sum on the left hand side we have, according to (1.2),

$$\frac{1}{1-t} \sum_{n=0}^{\infty} \left(\frac{t}{1-t}\right)^n \binom{\alpha}{n} H_n = \sum_{n=0}^{\infty} t^n \left\{ \sum_{k=0}^{n} \binom{n}{k}\binom{\alpha}{k} H_k \right\}. \tag{4.6}$$

Substituting this in (4.5) yields the representation (4.1). The proof is completed.

**Corollary 5.** For every integer $n \geq 0$ and every complex number $\alpha$ with Re $\alpha > -1$ we have

$$\sum_{k=0}^{n} \binom{n}{k}\binom{\alpha}{k} H_k = \binom{\alpha+n}{n}(H_\alpha + H_n - H_{\alpha+n}). \tag{4.7}$$

*Proof.* The identity follows immediately by comparing the representations (4.1) and

$$\frac{-\log(1-t)}{(1-t)^{\alpha+1}} = \sum_{n=0}^{\infty} (H_{\alpha+n} - H_\alpha) \binom{\alpha+n}{n} t^n, \tag{4.8}$$

that is, (2.15).

Note that the digamma function has the property

$$\psi(z+m) - \psi(z) = \frac{1}{z} + \frac{1}{z+1} + \ldots + \frac{1}{z+m-1} \tag{4.9}$$

and therefore,

$$H_\alpha + H_n - H_{\alpha+n} = \sum_{k=1}^{n}\left(\frac{1}{k} - \frac{1}{k+\alpha}\right) = \alpha \sum_{k=1}^{n} \frac{1}{k(k+\alpha)} . \tag{4.10}$$

## 5. More harmonic number identities

**Corollary 6.** For every complex $\alpha$ and every positive integer $n$,

$$\sum_{k=0}^{n} \binom{\alpha+1}{n-k}\binom{\alpha+k}{k}(-1)^{n-k} H_k = \frac{1}{n} + \binom{\alpha}{n}\frac{(-1)^{n-1}}{n} , \tag{5.1}$$

$$\sum_{k=0}^{n} \binom{\alpha}{n-k}\binom{\alpha+k}{k}(-1)^{n-k} H_k = H_n + \sum_{k=1}^{n}\binom{\alpha}{k}\frac{(-1)^{k-1}}{k} . \tag{5.2}$$

In particular, when $\alpha = n$, these identities become correspondingly,

$$\sum_{k=0}^{n} \binom{n+1}{n-k}\binom{n+k}{k}(-1)^{n-k} H_k = \frac{1+(-1)^{n-1}}{n} , \tag{5.3}$$

$$\sum_{k=0}^{n} \binom{n}{k}\binom{n+k}{k}(-1)^{n-k} H_k = 2H_n . \tag{5.4}$$

The identities (5.3) and (5.4) are very similar to the identities (3.123) and (3.122) in [6].

For the proof we need a simple lemma.

**Lemma 2.** For every power series as in (1,1) we have

$$\frac{1}{1-t}\sum_{n=0}^{\infty} a_n t^n = \sum_{n=0}^{\infty} t^n \left\{\sum_{k=0}^{n} a_k\right\} , \tag{5.5}$$

$$(1+\lambda t)^\alpha \sum_{n=0}^{\infty} a_n t^n = \sum_{n=0}^{\infty} t^n \left\{\sum_{k=0}^{n} \binom{\alpha}{n-k} a_k \lambda^{n-k}\right\} . \tag{5.6}$$

The proof is left to the reader. After expanding $(1+\lambda t)^\alpha$ one can either change the order of summation or use the Cauchy rule for multiplication of power series.

*Proof of the corollary.* We put $a_k = (-1)^{k+1} H_k$ in (2.4) to obtain

$$\sum_{n=0}^{\infty} \binom{\alpha}{n} H_n z^n = (z+1)^\alpha \sum_{n=0}^{\infty}\left(\frac{z}{z+1}\right)^n \binom{\alpha}{n}\frac{(-1)^n}{n} , \tag{5.7}$$

in view of the well-known binomial transform

$$\sum_{k=1}^{n} \binom{n}{k}(-1)^{k-1} H_k = \frac{1}{n} . \tag{5.8}$$

At the same time, from (4.3),

$$\sum_{n=0}^{\infty} \binom{\alpha}{n} H_n z^n + (1+z)^{\alpha} \log(1+z) = \frac{1}{1+z} \sum_{n=0}^{\infty} \left(\frac{z}{z+1}\right)^n H_n \binom{\alpha+n}{n}, \tag{5.9}$$

so that replacing the first sum here by the right hand side of (5.7) we obtain

$$(z+1)^{\alpha} \sum_{n=0}^{\infty} \left(\frac{z}{z+1}\right)^n \binom{\alpha}{n} \frac{(-1)^n}{n} + (1+z)^{\alpha} \log(1+z)$$

$$= \frac{1}{1+z} \sum_{n=0}^{\infty} \left(\frac{z}{z+1}\right)^n H_n \binom{\alpha+n}{n}. \tag{5.10}$$

With the substitution (1.4) this becomes

$$\sum_{n=0}^{\infty} t^n \binom{\alpha}{n} \frac{(-1)^{n-1}}{n} - \log(1-t) = (1-t)^{\alpha+1} \sum_{n=0}^{\infty} t^n H_n \binom{\alpha+n}{n}, \tag{5.11}$$

and according to (5.6), this can be written as

$$\sum_{n=0}^{\infty} t^n \binom{\alpha}{n} \frac{(-1)^{n-1}}{n} - \log(1-t) = \sum_{n=0}^{\infty} t^n \left\{ \sum_{k=0}^{n} \binom{\alpha+1}{n-k} \binom{\alpha+k}{k} (-1)^{n-k} H_k \right\}. \tag{5.12}$$

Using the expansion of the logarithm

$$-\log(1-t) = \sum_{n=1}^{\infty} \frac{t^n}{n} , \tag{5.13}$$

and comparing coefficients in (5.12) we arrive at (5.1). In order to obtain (5.2) we write (5.11) in an equivalent form, dividing by $1-t$,

$$\frac{1}{1-t} \sum_{n=0}^{\infty} t^n \binom{\alpha}{n} \frac{(-1)^{n-1}}{n} - \frac{\log(1-t)}{1-t} = (1-t)^{\alpha} \sum_{n=0}^{\infty} t^n H_n \binom{\alpha+n}{n}. \tag{5.14}$$

According to Lemma 3 and (3.10) this becomes

$$\sum_{n=0}^{\infty} t^n \left\{ \sum_{k=1}^{n} \binom{\alpha}{k} \frac{(-1)^{k-1}}{k} \right\} + \sum_{n=1}^{\infty} H_n t^n = \sum_{n=0}^{\infty} t^n \left\{ \sum_{k=0}^{n} \binom{\alpha}{n-k} \binom{\alpha+k}{k} (-1)^{n-k} H_k \right\}, \tag{5.15}$$

and yields (5.2). The proof is completed.